\documentclass[a4paper,12pt]{amsart}
\usepackage{epsfig,graphicx,verbatim}
\newcommand\urladdrx[1]{{\urladdr{\def~{{\tiny$\sim$}}#1}}}
\usepackage{natbib}
\bibpunct{[}{]}{,}{n}{,}{,}

\numberwithin{equation}{section}

\long\def\quote#1\endquote{{{\it#1}}}

\begin{document}
\title[David George Kendall, a biographical account]
{David George Kendall\\a biographical account}

\author{Geoffrey Grimmett}
\address{Statistical Laboratory,
Centre for
Mathematical Sciences,
University of Cambridge,
Wilberforce Road, Cambridge CB3 0WB, U.K.}
\email{g.r.grimmett@statslab.cam.ac.uk}
\urladdrx{http://www.statslab.cam.ac.uk/~grg/}

\date{October 6, 2008}

\begin{abstract}
This biographical account of the life and work of David Kendall
includes details of his personal and professional activities. 
Kendall is probably
best known for his work in applied probability, especially queueing theory,
and in stochastic analysis and spatial statistics.
\end{abstract}


\maketitle

\bigskip

David Kendall was a leading figure in probability and statistics,  
and a pioneer in applied probability, stochastic analysis 
and geometry, and associated statistical analysis.  Through his own 
research and that of his students in Oxford and Cambridge, he established 
modern probability theory in the United Kingdom.  With his recent death, 
we have lost another of the great generation of statisticians 
that spanned the 20th century.  

Kendall's skills and interests were broad.  His intention as a mathematics 
student was to become an astronomer.  The war intervened, and, 
following demobilization, he moved into the developing field of 
probabilistic modelling and analysis.  Making use of exceptional 
mathematical skills, he identified and solved problems in many areas 
including population modelling, branching processes, queueing theory, 
Markov processes, stochastic analysis, as well as stochastic geometry and 
the statistical analysis of geometrical data.

He undertook significant positions of leadership within the statistics and 
mathematics community, including being the first Professor of 
Mathematical Statistics at Cambridge University from 1962--1985, and 
the President of the Bernoulli Society for Mathematical Statistics and 
Probability at its inauguration in 1975.

\begin{figure}[ht]
\begin{center}
\includegraphics[angle=0,width=8cm]{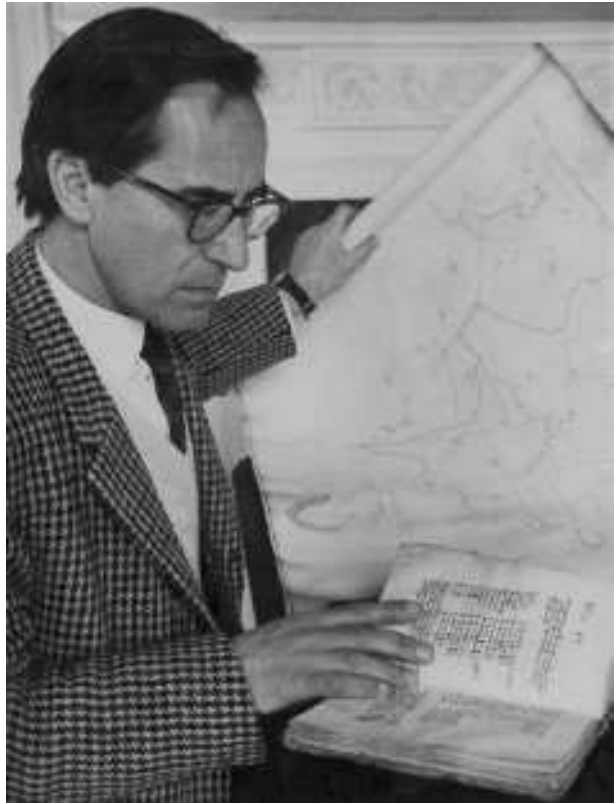}
\caption{David George Kendall.
Probably taken around 1975 when DGK was working on contiguity data
from manorial records of Whixley, Yorkshire.}
\label{DGk-photo}
\end{center}
\end{figure}

David Kendall was born in Ripon on 15 January 1918, 
the son of Fritz Ernest Kendall and
Emmie Taylor. He remained faithful to his 
Yorkshire roots, spending frequent holidays there in later life.  
He was influenced as a teenager by 
George Viccars, the senior mathematics master at Ripon Grammar 
School. His interest in astronomy was sparked in the early 1930s by a 
series of radio talks by James Jeans.  He realised that a strong 
mathematical background was essential to becoming a successful 
astronomer, and he applied to Oxford and Cambridge accordingly. 
Gonville and Caius College's offer of an Exhibition could not be 
accepted on financial grounds, and instead he accepted in 1936 a (closed) 
Hastings Scholarship at Queen's College, Oxford, where his tutor was
U.\ S.\ Haslam-Jones.  He formed ties to H.\ H.\ Plaskett and E.\ A.\ Milne, 
and, under their influence, he published his first research paper at the end 
of his second undergraduate year.  It dealt with the computation of a 
certain integral, and it appeared 
in the {\it Zeitschrift f\"ur Astrophysik\/}, \cite{DGK0}.

On graduating in 1939 with a first class degree, Kendall was appointed to 
a Senior Studentship in Astronomy at Balliol College.  The outbreak of 
war disturbed his plans, and he left his new post in 1940 for work on 
rocketry at the Projectile Development Establishment (PDE) in West 
Wales.  When in due course Frank Anscombe and Maurice Bartlett 
moved away from PDE, Kendall's then boss Louis Rosenhead allowed him one week's 
leave to learn statistics.  With the help of private tuition from Anscombe, a week 
proved sufficient for the matters at hand.   
Kendall has described (with Kenneth Post) his work on guidance systems 
for rockets in \citep{KP96,KP97}.
After 
the declaration of peace, he travelled with T-force to talk to German 
scientists, thereby procuring in G\"ottingen a photograph of Gauss on his 
death bed.  

Kendall returned to Oxford in 1946 as Tutorial Fellow of Magdalen 
College, a post he was to hold until his move to Cambridge in 1962.  He 
was heavily influenced during this period by Maurice Bartlett.  Amongst 
his early works of significance are two papers presented at discussion 
meetings of the Royal Statistical Society (RSS), 
dealing respectively with models for the growth of 
populations, and queues.  He spoke with Moyal and Bartlett at the RSS
Symposium on Stochastic Processes in 1949, and contributed a long and 
significant paper (with discussion) on stochastic models for population 
growth, \cite{RSS1}.  
This partly expository paper includes accounts of his own 
investigations, and it concerns one of his principal interests of the 
next decade, 
namely the theory of continuous-time Markov processes.  His second paper 
\cite{RSS2} summarises the state of queueing theory in 1951.  Dennis Lindley's 
contribution to the discussion is especially interesting: firstly, he 
proposed the analysis of the waiting-times of successive customers in 
order to solve the general G/G/1 queue, but without explicit mention of 
the duality relation leading to the ladder-height representation; and secondly, 
he took up the notion of a `regeneration point', to which Kendall and his 
collaborators were to return later.  

One of David's useful contributions to queueing theory was the D/G/$k$ 
notation of his 1953 paper in the {\it Annals of Mathematical Statistics\/}, 
\cite{Ken53},
and 
he would jest about his Guy Medal in Silver, awarded 
by the RSS in 1955, as 
\lq\lq perhaps the only medal to be awarded for inventing a terminology and a 
notation\rq\rq.  (He reminisced freely about his life in conversation with Nick 
Bingham, see \cite{Bingh}.) His major lasting 
contribution to queueing theory was probably the analysis via embedded 
Markov chains.  It was around 1953 that David met Harry Reuter, thus 
initiating a happy collaboration into the foundational aspects of 
continuous-time Markov chains.  Quite a few of David's interests from 
the 1950s are well represented in current undergraduate and Masters-level 
lecture courses around the world.  
Daryl Daley and David Vere-Jones have written about
David's contributions to applied probability in \cite{Bin}.

Kendall was elected in 1962 to the new Professorship of Mathematical Statistics 
at Cambridge, and to a Professorial Fellowship at Churchill College.  
There was evidently some disappointment in 
certain quarters that someone more obviously a statistician had not been 
found, but David's candidacy was overwhelmingly strong.  He obliged by 
working on several problems involving analysis of data, while continuing 
his more theoretical work.  His election to Cambridge permitted a 
regeneration of the Statistical Laboratory, and even led to a courtesy visit 
by Ronald Fisher to its then home in the basement of the Chemistry Building.
He remained the {\it ex officio\/} Director of the Lab until the strings were loosened
in 1973 and Peter Whittle took over leadership. 

David's choice of problem area was characteristically individualistic, 
including inference problems of archaeology, such as grave sequencing 
and the reconstruction of local maps from contiguity data, and the 
distribution of standing stones and the statistics of ley lines.  He 
continued working on models of applied probability inspired by real-life 
problems.  Two examples: a Brownian motion model to reconstruct the optic nerve 
connections following trauma; and modelling the flight of a bird equipped 
with \lq\lq a clock, $\dots$ a sextant, $\dots$ and a copy of what used to be called the 
{\it Nautical Almanac\/}\rq\rq.  At the same time, he published several works of a 
more analytical nature: on `Delphic'  semigroups, and on regeneration 
and renewal.  

In much of his work he was able to combine stochastic 
analysis and stochastic geometry, two strands that will forever be 
associated with his name through the two volumes \cite{KH1,KH2}
 in commemoration of 
Rollo Davidson, edited by DGK and E.\ F.\ Harding and published in 1973--1974.
David was greatly saddened by the loss of Rollo in an Alpine climbing 
accident in 1970 and
he was instrumental with colleagues in the establishment of the Rollo Davidson Trust at Churchill
in 1975. 
From Rollo's obituary in {\it Stochastic Analysis\/}/{\it Geometry\/}:
\lq\lq The hazards of [mountain adventure], never wholly to be avoided, are familiar to all, 
and to rail at its folly is to invite a reply which he himself [Rollo] might have made, 
$\dots$, 
`If you always look over your shoulder, how can you still remain a human being?'.\rq\rq 

In later years spanning his retirement in 1985, David worked on a theory 
of shape that extended his work on ley lines.  This culminated in 1995 in 
the publication of the book {\it Shape and Shape Theory\/}, 
\cite{Ketal}, written 
by DGK, D.\ Barden, T.\ K.\ Carne, and H.\ Le, and described by a reviewer for 
Mathematical Reviews as \lq\lq a mathematical jewel\rq\rq.

David Kendall was much decorated, receiving amongst other awards the 
Guy Medal in Silver (1955) and in Gold (1981) of the Royal Statistical Society, 
the Sylvester Medal 
(1976) of the Royal Society, and both the Senior Whitehead Prize (1980) 
and the De Morgan Medal (1989) of the London Mathematical Society.  
He was elected to the Fellowship of the Royal Society (1964), to 
Membership of the Academia Europaea (1991), to Honorary 
Membership of the Romanian Academy (1992), and to the Fellowship
of the Institute of Mathematical Statistics.  He received honorary 
degrees from the Universit\'e Ren\'e Descartes (1976), and the University of 
Bath (1986), and he took his Oxford DSc in 1977.

He was a mainstay of several learned societies, including service as the 
first President of the Bernoulli Society in 1975, and as President of the 
London Mathematical Society (LMS) from 1972--1974. He was a Council member
of the Royal Statistical Society from 1950--1952, and served as
Vice-President from 1963--1964. 
When Joe Gani, Ted Hannan, and Norma McArthur
were striving for support to set up what became 
the Applied Probability Trust in 1963, it
was David who provided the key introduction to 
the LMS.
He became the first Professorial Fellow elected 
by Churchill College, Cambridge, in 1962, and he remained a Fellow 
until his death.  He was elected to an Honorary Fellowship of 
his old Oxford college, Queen's, in 1985.

One might describe David as a gentleman of the old school.  
On occasion he could seem distant, but, on 
penetrating his slightly formal manners, one found a person brimming 
with knowledge and
human values, and eager to help younger people with common interests.
A devout Anglican throughout his life, David attended services 
at the Church of St Mary and St Michael in Trumpington.  Through his 
students his scientific legacy, already substantial, has been multiplied 
manyfold.  It may seem invidious to select a few from a list of the august, 
but, as an indication of the scientific breadth of his influence, we mention 
David Edwards, David Williams, Daryl Daley, David Vere-Jones,
John Kingman, Rollo Davidson, Nick Bingham, Richard Tweedie, 
and Bernard Silverman.

\begin{figure}[ht]
\begin{center}
\includegraphics[angle=0,width=11cm]{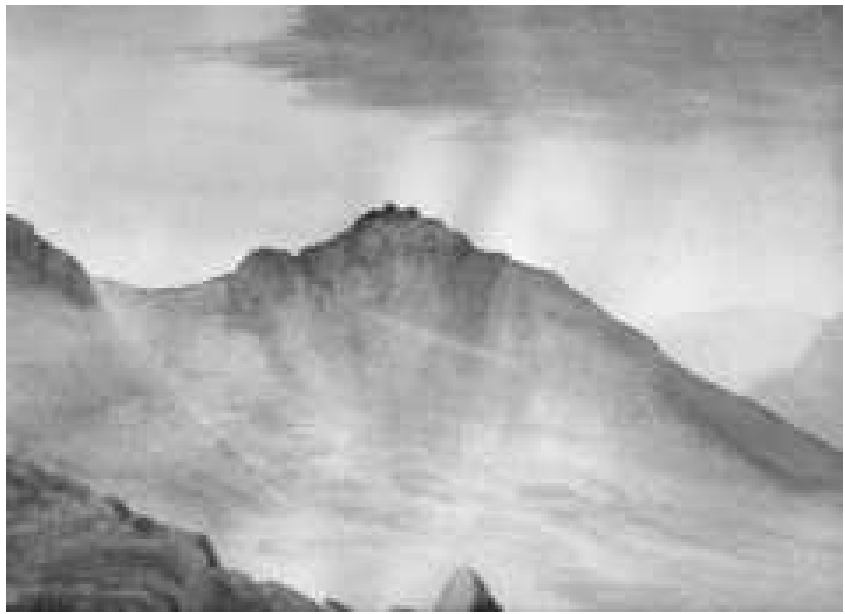}
\caption{The east face of Tryfan, by DGK in watercolour. Adam and Eve may be
distinguished, and also the Heather Terrace. Reproduced
by kind permission of the Kendall family.}
\label{Tryfan}
\end{center}
\end{figure}

David travelled widely, forming personal connections that were to be 
maintained.  In later years he tended eastwards, making many trips to the 
countries of Eastern Europe and further afield.  From his famous visit to 
China in 1983 with David Williams, he returned with a student, Huiling Le, 
now at Nottingham University.  Many will recall 
his effective captaincy of the British team at that most extraordinary First 
World Congress of the Bernoulli Society held in Tashkent in 1986.  He 
had a great love of the outdoors, and he took much pleasure from sharing 
climbing yarns and examining fossils, collected preferably from the cliffs at Whitby.  
One might easily believe that he 
accorded greater credit to the Reverend Henry William Watson for his 
founding membership of the Alpine Club than for his (incomplete) 
solution to the extinction problem for branching processes. 

David did not permit retirement to intervene in his intellectual activities.  
He remained for many years a presence in the Stats Lab, until the distance 
to the new site in Clarkson Road became a hindrance and he gave up 
cycling.  Until shortly before his death he could be seen striding 
purposefully around Cambridge, and he frequently attended Lab parties 
and College lunches.  It was with sadness that his colleagues and friends 
learned of his death on 23 October 2007, following a chronic illness 
through which he was nursed by his wife Diana.  The University Church 
of Great Saint Mary's was comfortably filled by the celebrants at his 
Memorial Service, and the crowd convened to Churchill for tea and
tributes. 

David and Diana (n\'ee Diana Louise Fletcher) married in 1952. 
They have six children, 
Wilfrid, Bridget, Felicity, Judy, and twins George and Harriet, spanning 
just seven years in age difference, and  
eight grandchildren at the time of 
writing.

\newcommand\vol{\textbf}
\newcommand\jour{\emph}
\newcommand\book{\emph}
\newcommand\inbook{\emph}
\def\no#1#2,{\unskip#2, no. #1,} 
\newcommand\toappear{\unskip, to appear}

\end{document}